\documentclass{article}
\usepackage[T1]{fontenc}
\usepackage[latin1]{inputenc}
\usepackage{latexsym}                 
\usepackage{amssymb}
\usepackage{amsthm}
\usepackage{color}
\usepackage{graphicx}
\usepackage{hyperref}
\usepackage{verbatim}
\def\eref#1{(\ref{#1}%
)}
\def\RSref#1{\ref{#1}%
}
\def\RSlabel#1{\label{#1}%
}
\def\RScite#1{\cite{#1}%
}

\newcommand{\bql}[1]{%
\begin{equation}\label{#1}%
}
\def\filename#1{}
\newcommand{\eq}{\end{equation}}
\def\fa{\hbox{ for all }}

\def\dfrac#1#2{\displaystyle{\frac{#1}{#2}   }}

\newcommand{\R}{\ensuremath{\mathbb{R}}}

\newcommand{\red}[1]{{\color{red}{#1}}}    
 
\def\calm{\mathcal{M}}

\def\biglf{\par\bigskip\noindent}
\def\bc{\mathbf c}
\def\bbf{\mathbf f}
\def\bg{\mathbf g}
\def\br{\mathbf r}
\def\bs{\mathbf s}
\def\bu{\mathbf u}
\def\bA{\mathbf A}
\def\bB{\mathbf B}
\def\bC{\mathbf C}
\def\bR{\mathbf R}
\def\bS{\mathbf S}
\def\bU{\mathbf U}
\def\bbf{\mathbf f}
\begin{document}
\begin{center}
{\Large \bf A computational tool for comparing}
\par\medskip
{\Large \bf all linear PDE solvers}
\biglf
{\large \bf - Optimal methods are meshless -}
\biglf
Robert Schaback \\
Univ. G\"ottingen\\
\biglf
Draft of \today\\
\biglf
schaback@math.uni-goettingen.de
http://num.math.uni-goettingen.de/schaback/research/group.html
\end{center}
\biglf 
{\bf Abstract}: The paper starts out with a 
computational technique that allows to compare
all linear methods for PDE solving that use the same input data.
This is done by writing them as linear recovery formulas for  solution values
as linear combinations of the input data. Calculating the norm
of these reproduction formulas on a fixed Sobolev space will then 
serve as a quality criterion that allows a fair comparison
of all linear methods with the same inputs, 
including finite--element, finite--difference and meshless local
Petrov--Galerkin techniques.  A number of illustrative examples will be
provided. As a byproduct, it turns out that a unique error--optimal method exists.
It necessarily outperforms 
any other competing technique using the same data, e.g. those just mentioned,
and it is necessarily meshless, if solutions are written 
``entirely in terms of nodes'' (Belytschko et. al. 1996 \RScite{belytschko-et-al:1996-1}).
On closer inspection, it turns out that it 
coincides with {\em symmetric meshless collocation} carried 
out with the kernel of the Hilbert space used for error evaluation, e.g. with
the kernel of the Sobolev space used. This technique 
is around since at least 1998, but its optimality properties went unnoticed, so far.
\section{Introduction}\RSlabel{SecIntro}
For simplicity, consider a standard Dirichlet problem
\bql{eqLufBug0}
\begin{array}{rcll}
Lu&=&f & \hbox{ in } \Omega\\ 
Bu&=&g & \hbox{ in } \Gamma:=\partial \Omega\\ 
\end{array}
\eq
with a linear differential operator $L$ and a linear boundary operator $B$
posed on a domain $\Omega$. Assume that the amount of available information is
fixed and limited, namely to values 
\bql{eqdata}
\begin{array}{rlll}
f(x_j) & \hbox{ in points } & x_j\in \overline \Omega,& 1\leq j\leq m,\\
g(y_k) & \hbox{ in points } & y_k\in \Gamma:=\partial\Omega,& 1\leq k\leq n.\\
\end{array} 
\eq
Under all methods for solving such problems, we ask for the one with smallest
error that uses this information and not more. To make this more precise,
we fix a single point $x\in\overline\Omega$ and assume existence of a true
solution $u^*$ of the problem. Any method $\calm$ using the above information will
produce some numerical solution $\tilde u_\calm$, and we want to single out methods that 
make the error $|u^*(x)-\tilde u_\calm (x)|$ small for all 
problems posed that way. We achieve this by looking at error bounds 
$$
|u^*(x)-\tilde u_\calm(x)|\leq C _\calm\|u^*\|_H  \fa u^* \in H,
$$
where the problems and their solutions
are allowed to vary in such a way that the true solutions lie in some fixed 
reproducing kernel Hilbert
space $H$ of functions on $\Omega$, e.g. a fixed Sobolev space. 
In this sense, the constant $C_\calm$ 
describes the worst--case error behavior of the method $\calm$ on all problems
with solutions in $H$ using the same data. We shall show how to evaluate $C_\calm$ 
numerically, and this will allow fair comparisons between methods using the same
data. A few examples including important methods like finite elements,
gerneralized finite differences,  and 
meshless collocation in various forms are provided at the end.
\biglf
But one can also ask for a method $\calm$ that 
makes $C_\calm$ minimal over all linear methods using the same data. 
This problem will be written as one of optimal recovery of
functions, and it is proven that 
the optimal solution exists uniquely. 
It can be calculated explicitly, takes the form of a recovery method in
reproducing kernel Hilbert spaces, and is a {\em meshless method}, if users ask for 
values of the solution at nodes.
Since it is optimal in the above sense, it outperforms errorwise
any other competing technique using the same data, may it use finite elements,
finite differences, or local Petrov--Galerkin techniques. Upon
closer inspection, it turns out to be nothing new, since it 
is a special case of {\em symmetric meshless collocation},
computed using the kernel that was used for error assessment.
This method relies on Hermite--Birkhoff interpolation 
in reproducing kernel Hilbert spaces \RScite{wu:1993-2} and its application
to PDE solving was 
first analyzed in \RScite{franke-schaback:1998-2a,franke-schaback:1998-1}.
\section{Recovery Problems}\RSlabel{SecRP}
To keep the presentation simple, we stay with the Dirichlet problem example.
Any linear method that provides an approximate solution value $\tilde u(x)$ for
a fixed point $x$ and uses exclusively the data \eref{eqdata}, must necessarily 
satisfy a formula of the form
\bql{eqursg}
\tilde u(x)=\displaystyle{\sum_{j=1}^m r_j(x)f(x_j) +\sum_{k=1}^ns_k(x)g(y_k),   }  
\eq
whatever the weights $r_j(x)$ and $s_k(x)$ are.  We call this
a {\em direct linear recovery}  of the value $\tilde u(x)$ from the data
on the right--hand side. 
The existence of \eref{eqursg} follows from linearity and the restriction to the admitted data,
but in general it will need some effort to
rewrite a classical method in this form. We come back to this in the 
section on specific methods.
\biglf
Clearly, the error in \eref{eqursg} is
$$
u^*(x)-\tilde u(x)
=u^*(x)-\displaystyle{ \sum_{j=1}^mr_j(x)\Delta u^*(x_j) -\sum_{k=1}^ns_k(x)u^*(y_k)  }, 
$$
where we now have connected the data functions $f$ and $g$ back to the true
solution $u^*$. The map
\bql{eqexrs}
\epsilon_{x,\br,\bs}\;:\; u^*\mapsto 
u^*(x)-\displaystyle{ \sum_{j=1}^mr_j(x)\Delta u^*(x_j) -\sum_{k=1}^ns_k(x)u^*(y_k)  }
\eq
is a linear functional in $u^*$, and if it is bounded on some Hilbert space $H$,
we have an error bound
\bql{equueHuH} 
|u^*(x)-\tilde u(x)|\leq \|\epsilon_{x,\br,\bs}\|_{H^*}\|u^*\|_H \fa u^*\in H.
\eq
If the norm $\|\epsilon_{x,\br,\bs}\|_{H^*}$ is evaluated, it precisely
describes the worst--case
error behavior for all problems with solutions in $H$, 
because the inequality is sharp by definition of the norm of the functional.
In other words, one can express the error explicitly in percent of $\|u^*\|_H$,
and this is what we shall do in the rest of the paper. 
Furthermore, one can ask for a choice $\br^*$ and $\bs^*$ of the coefficients that
make the norm minimal, and this will later lead to an optimal method for the given
reconstruction. 
\biglf
Note that we do not consider the way numerical solutions of the considered
methods are actually calculated.
The amount of numerical integration needed for weak methods, 
the error committed in integration
subroutines, and any multigrid solver techniques are completely irrelevant.
Different numerical strategies within the same class of methods, e.g. finite element
solvers with
different ways of integration or with different element families can be
compared directly as well. The only requirement is that the method is rewritten
in direct linear recovery form \eref{eqursg}, and then the total error is evaluated, containing
all method--specific internal features of the technique. 
\section{Error Evaluation}\RSlabel{SecEE}
Fortunately, it is easy to evaluate such norms in reproducing kernel
Hilbert spaces, and Sobolev spaces are examples. In such spaces, the inner
product and the
kernel $K\;:\;\Omega\to\Omega$ have the properties
$$
\begin{array}{rcl}
f(x)&=&(f,K(x,\cdot))_H \fa f\in H, \;x\in\Omega\\
K(x,y)&=&(K(x,\cdot),K(y,\cdot))_H\fa x,y\in\Omega\\
\end{array}
$$
and for all linear and continuous functionals $\lambda,\mu\in H^*$, 
$$
\begin{array}{rcl}
(\lambda,\mu)_{H^*} &=& (\lambda^xK(x,\cdot),\mu^yK(y,\cdot))_H\\
&=& \lambda^x\mu^yK(x,y),
\end{array}
$$ 
where the upper index at the functionals denotes the variable that is acted
upon. Details are in the literature on kernels, with \RScite{wendland:2005-1}
being a fairly complete reference and \RScite{schaback:1997-3,schaback:2011-1} being
open access compilations for teaching purposes.
\biglf
Global Sobolev spaces $W_2^m(\R^d)$ in $d$
dimensions
and for order $m>d/2$ in the standard Fourier transform definition
have the radial {\em Mat\'ern} reproducing kernel
\bql{eqKnu}
K(x,y)=\dfrac{2^{1-m}}{\Gamma(m)}\|x-y\|_2^{m-d/2}K_{m-d/2}(\|x-y\|_2),\;x,y\in\R^d
\eq
with the modified Bessel function $K_{m-d/2}$ of the second kind. Local Sobolev
spaces $W_2^m(\Omega)$ for domains $\Omega\subset\R^d$  are norm--equivalent
to the global spaces, as long as domains are non--pathological, i.e. they 
satisfy a Whitney extension property or a have a piecewise smooth boundary with
a uniform interior cone condition. We shall use the above kernels for evaluation
of errors in Sobolev space.
\biglf
To avoid double sums and to pave the way for generalizations, we  
shall calculate the error norms after 
rewriting \eref{eqexrs} in terms of functionals as
$$
\epsilon_{x,\bc}=\delta_x-
\displaystyle{ \sum_{k=1}^{N}c_k(x)\lambda_k  }
$$
with
\bql{eqlamlam}
\lambda_j(u)=L u(x_j),\;1\leq j\leq m, \hbox{ and } 
\lambda_{m+i}(u)=u(y_i),\;1\leq i\leq n,\;N=m+n
\eq 
in our special case. It should be clear that other differential operators
and other boundary conditions will
just change the functionals here.
\biglf
Then we get the quadratic form
\bql{eqexcsquare}
\begin{array}{rcl}
\|\epsilon_{x,\bc}\|^2_{H^*}
&=&
K(x,x)-2\displaystyle{ \sum_{i=1}^{N}c_i(x)\lambda_i^zK(x,z)}\\
&&+
\displaystyle{\sum_{j,i=1}^{N}c_i(x)c_j(x)\lambda_i^y\lambda_j^zK(y,z)   } 
\end{array}
\eq
which can be explicitly evaluated if the kernel is known
and the functionals are continuous, in particular on
all Sobolev spaces on which the functionals are continuous. For instance,
if $L$ is a second--order elliptic operator, this poses the restriction $m-2>d/2$ on the 
Sobolev order $m$ we can use. But this is normal, since we focus on methods
that use $f=Lu^*$ pointwise, forcing $f$ to be in $W_2^k$ with $k>d/2$ and thus
$u^*\in W_2^m$ with $m>2+d/2$. All methods that use these data
are implicitly making this smoothness assumption, even if their users think that they
are working in in less regular spaces. Plenty of papers and books miss this
point.
In particular, the standard FEM method is usually formulated 
and analyzed in low--regularity spaces
like $H^1$ or $H^2$, but when it comes to implementations
using data like \eref{eqdata}, it needs $H^m$ with $m>3$ in $\R^2$.
\biglf
All methods based on the data $\lambda_1(u^*),\ldots,\lambda_K(u^*)$
and brought into the recovery form
$$
\tilde u(x)=\displaystyle{\sum_{i=1}^{N}c_i(x)\lambda_i(u^*)   } 
$$
generalizing \eref{eqdata} and  
\eref{eqexrs} can now be plugged into 
\eref{eqexcsquare} to show how good the reproduction quality at $x$ is,
since \eref{equueHuH} generalizes accordingly. Note that there is no linear
system to be solved once the recovery formula is known. But since the 
formula is inserted into the positive definite quadratic form
\eref{eqexcsquare}, a
small final value will necessarily contain quite some amount of cancellation.
The computational complexity for norm evaluation is ${\cal O}(N^2)$.
\biglf
This allows a fair comparison of all such methods using the same data, 
and we shall provide examples
below. The comparison can be made pointwise, as we saw, but 
for small problems one can plot the function 
$x\mapsto \|\epsilon_{x,\bc}\|^2_{H^*}$
to see where a method works badly and needs more data. This gives direct
information  for refinement of the discretization.
\section{Optimal Methods }\RSlabel{SecOM}
Since \eref{eqexcsquare} is a quadratic form with a positive semidefinite
matrix, it can be minimized. The optimal solution coefficients $c_i^*(x)$ solve
the system
\bql{eqcllKlK}
\displaystyle{\sum_{k=1}^Nc_k^*(x)\lambda_j^y\lambda_k^zK(y,z)=\lambda_j^z
  K(x,z),\;1\leq j\leq N,   } 
\eq
and the system can be proven to be solvable \RScite{wu:1993-2}. Clearly, this
choice of coefficients will then outperform all other competitors error--wise.
We shall show examples later. The minimal value of \eref{eqexcsquare} then is 
$$
\begin{array}{rcl}
\|\epsilon_{x,\bc^*}\|^2_{H^*}
&=&
K(x,x)-\displaystyle{ \sum_{i=1}^{N}c_i(x)^*\lambda_i^zK(x,z)}\geq 0
\end{array}
$$
and does not contain any matrix.
\biglf
The system \eref{eqcllKlK} reveals the nature of this method. 
Indeed, the functions $c_k$ are necessarily linear combinations of the functions
$\lambda_i^zK(\cdot,z)$, and application of $\lambda_i^x$ to the above system 
shows the Lagrange property $\lambda_i(c_k)=\delta_{ik},\;1\leq i,k\leq N$. This means
that the $c_k$ are the  Lagrange basis for general Hermite--Birkhoff interpolation
of the given data by the functions $\lambda_i^zK(\cdot,z)$,,and this is 
the well--known method of {\em symmetric meshless collocation} based
on \RScite{wu:1993-2} and analyzed thoroughly in 
\RScite{franke-schaback:1998-1,franke-schaback:1998-2a}. 
The optimality of the technique in the sense of this paper 
should have been known since at least since 1997, but it went unnoticed because
\cite[p.82, (4.2.2)]{schaback:1997-3} was not applied to PDE solving
at that time.
\section{Special Methods}\RSlabel{SecSM}
We now give some details on how to evaluate recovery errors for special PDE
solution techniques. This will be useful for the examples in the final section,
and we specialize to $L=-\Delta $ here.
\subsection{Finite Elements}\RSlabel{SecSMFEM}
The simplest possible 2D finite--element code uses piecewise linear elements on the
triangles of a triangulation of a domain with piecewise linear boundary,
and requires $f$--values only at the barycenters of the triangles. These are the
points $x_j$ in \eref{eqdata} in the FEM version that we denote by FEMBary below.  But since
usually there are more triangles than nodes, one can also prescribe $f$--values
at all interior and boundary nodes to calculate approximate values at the
barycenters. This usually needs less $f$--values at the same order of accuracy. We call
this method FEMNode below. These two FEM variations have different recovery
formulas \eref{eqursg} and different error functionals \eref{eqexrs} to be
compared. 
\biglf
To get the recovery formulas in the form needed here, users will have to check carefully what
their FEM code does. We used the MATLAB {\tt pdetool} setting, which does the
following. It uses all $N$ triangle vertices $z_1,\ldots,z_N$ for setting up 
the standard piecewiese linear test and trial functions $v_1.\ldots,v_N$ 
with $v_j(z_k)=\delta_{jk}$ and builds
the $N\times N$ stiffness matrix with entries $(\nabla v_i,\nabla v_j)_{L_2}$
in the usual way. The right--hand sides $(f,v_j)_{L_2},\;1\leq j\leq N$ are
certain linear combinations of $f$--values either at triangle barycenters or at
nodes of the neighboring triangles. 
These linear combinations form a sparse {\em integration matrix} $\bB$ with
entries
$b_{jk}$ such that the stiffness system without boundary conditions 
but with some form of numerical integration would be  
$$
\displaystyle{\sum_{i=1}^N (\nabla v_i,\nabla v_j)_{L_2}u(z_i) =
\sum_{k=1}^Nb_{jk}f(x_k),\;1\leq j\leq N.   } 
$$
The FEMBary and FEMNode variations have different $\bB$ matrices and use
$f$ at different points $x_k$, but the stiffness matrices are the same.
\biglf
The unknowns $u(z_i),\;i\in D\subset\{1,\ldots,N\}$ are known Dirichlet
values, and thus
$$
\displaystyle{\sum_{i\notin D} (\nabla v_i,\nabla v_j)_{L_2}u(z_i) =
\sum_{k=1}^Nb_{jk}f(x_k)-\sum_{i\in D} (\nabla v_i,\nabla v_j)_{L_2}u(z_i) ,\;j\notin D   } 
$$
is the system to be actually solved. Since the $u(z_i)$ with $i\in D$ are
$g$--values,
it has the matrix--vector form
\bql{eqAuBfCg} 
\bA\bu=\bB\bbf+\bC\bg
\eq
under adequate notation, and each row of 
$$
\bu=\bA^{-1}\bB\bbf+\bA^{-1}\bC\bg
$$
provides one instance of \eref{eqursg} for each of the points $x=x_i,\;i\notin
D$.
If we focus on the origin as one of these points, we need just one row.
\biglf
If users want the error at a non--nodal point $x$, they have to add a piecewise
linear interpolation on a triangle containing $x$, and then 
the discrete reconstruction formula is a linear combination of
three rows of the above system. This illustrates how numerical integration and
interpolatory post--processing both  enter explicitly into our version of
a complete and explicit FEM error analysis.
\subsection{Symmetric Kernel--Based Collocation}\RSlabel{SecSubSKBC}
These work on the two point sets  $X$ and $Y$ from
\eref{eqdata} and use linear combinations
\bql{equcDKdK}
u(x)=\displaystyle{\sum_{j=1}^mc_j\Delta K(x,x_j)+\sum_{k=1}^n d_kK(x,y_k)  } 
\eq
 of basis functions derived from a smooth kernel $K$. The argument in
Section \RSref{SecOM} shows that these methods yield optimal errors in the
Hilbert spaces in which their kernels are reproducing.  In case of
Sobolev spaces, we thus get the error--optimal methods this way.
\biglf
The numerical process collocates these trial functions at 
Dirichlet and PDE nodes, forming the block system
$$
\begin{array}{rcl}
\displaystyle{\sum_{j=1}^mc_j\Delta^x\Delta^y
  K(x_i,x_j)+\sum_{k=1}^n d_k\Delta K(x_i,y_k)}
&=& f(x_i),\;1\leq i\leq m,\\
\displaystyle{\sum_{j=1}^mc_j\Delta K(y_i,x_j)+\sum_{k=1}^n d_kK(y_i,y_k)}
&=& g(y_i),\;1\leq i\leq n.
\end{array}
$$
The inverse of the coefficient matrix recovers the coefficients $c_j$ and $d_k$
from the 
$f$ and $g$ data, and the numerical solution at $x$ is just a linear combination
\eref{equcDKdK} of those coefficients. 
Consequently, the discrete recovery \eref{eqursg} is furnished by
the
inverse of the above ``stiffness'' matrix, premultiplied by the row vector
of the kernel values in \eref{equcDKdK}. 
\biglf
Of course. one can use a special kernel $K$ to calculate the discrete recovery,
and then use
another kernel, e.g. one generating a Sobolev space, for error evaluation on
that space. We shall do this in the final section.
\subsection{Unsymmetric Kernel--Based Collocation}\RSlabel{SecSubUSKBC}
In contrast to the previous section,
this class of methods takes an additional set $Z=\{z_1,\ldots,z_N\}$ of usually $N=m+n$ points
and works on linear combinations
\bql{equdK}
u(x)=\displaystyle{\sum_{k=1}^N d_kK(x,z_k)  } 
\eq
of basis functions, while still using collocation in the point sets $X$ and $Y$
from \eref{eqdata}. This method dates back to early papers of Ed Kansa
\RScite{kansa:1986-1,kansa:1990-1} 
and was called MLSQ2 as a variation of the Meshless Local Petrov--Galerkin
method \RScite{atluri-zhu:1998-1,atluri:2005-1} of S.N. Atluri and collaborators. The linear system 
for the coefficients now is
$$
\begin{array}{rcl}
\displaystyle{\sum_{k=1}^n d_k\Delta K(x_i,z_k)}
&=& f(x_i),\;1\leq i\leq m,\\
\displaystyle{\sum_{k=1}^n d_kK(y_i,z_k)}
&=& g(y_i),\;1\leq i\leq n,
\end{array}
$$
and the inverse of the coefficient matrix (if it exists, see
\RScite{hon-schaback:2000-1}), premultiplied
by the vector of values $K(x,z_k)$ of \eref{equdK} will yield a row vector
for the discrete recovery formula \eref{eqursg}. If $N$ is chosen larger than
$m+n$ to increase stability, a pseudoinverse of the system coefficient matrix
can replace the inverse. This was done in the examples of the final section.
\subsection{Meshless Lagrange Methods}\RSlabel{SecMLM}
Here, a set $Z=\{z_1,\ldots,z_N\}$ of trial nodes is chosen, and there are
{\em shape functions} $u_1,\ldots,u_N$ such that {\em trial functions}
$$
u(x)=\displaystyle{\sum_{k=1}^Nu_k(x)u(z_k)   } 
$$ 
can be written ``entirely in terms of nodes''. Usually, this implies
Lagrange conditions $u_j(z_k)=\delta_{jk}$, and in many cases the shape functions
are
defined via Moving Least Squares. We do not care here for details,
and allow such 
techniques to come in weak or strong form. The strong case collocates in sets $X$
and $Y$ like above, forming a system
$$
\begin{array}{rcl}
\displaystyle{\sum_{k=1}^N\Delta u_k(x_j)u(z_k)   } &=& f(x_j),\;1\leq j\leq m,\\
\displaystyle{\sum_{k=1}^Nu_k(y_i)u(z_k)   } &=&  g(y_i),\;1\leq i\leq n \\
\end{array} 
$$
and the weights of the discrete recovery at $z_k$ will be a row of the
preudoinverse of the coefficient matrix of this system. 
\biglf
The weak cases form stiffness matrices and right--hand sides
like in the FEM situation, and then we get the coefficients of the
discrete recovery in the same way, involving a special matrix $\bB$ 
caring for the numerical integration.
\subsection{Generalized Finite--Difference Methods}\RSlabel{SecGDM}
Here, there are no trial functions, but everything is still expressed 
in terms of values at nodes $Z=\{z_1,\ldots,z_N\}$. In the strong situation,
PDE operator values are approximated by formulas like
\bql{eqDuau} 
\Delta u(x_j)\approx \displaystyle{\sum_{i=1}^N\alpha_{jk}u(z_k),   } 
\eq
with localized weights $\alpha_{jk}$, 
and this can be done with minimal error in a reproducing kernel Hilbert space
using  the logic of section \RSref{SecOM}. See 
\RScite{davydov-schaback:2012-1,schaback:2013-1} for more details. 
The linear system then is 
$$
\begin{array}{rcl}
f(x_j)
&=& \displaystyle{\sum_{i=1}^N\alpha_{jk}u(z_k),   }\\
f(x_j)-\displaystyle{\sum_{i\in D}\alpha_{jk}u(z_k), }
&=& \displaystyle{\sum_{i\notin D}\alpha_{jk}u(z_k),   }\\
\end{array} 
$$
if we use a subset $D$ of Dirichlet nodes like in the FEM case. This is of the
form \eref{eqAuBfCg} and we already know how to derive the
recovery formulas in such a case. 
It is interesting to see that the matrix with coefficients $\alpha_{jk}$ 
plays the role of a stiffness matrix here, and it is the place where sparsity 
can be implemented to yield local methods with sparse matrices. This was 
very successfully done 
in various application papers of C.S. Chen, B. Sarler, and G.M. Yao 
\RScite{sarler:2007-1,yao-et-al:2010-1, yao-et-al:2011-1,%
yao-et-al:2012-1,vertnik-sarler:2011-1}.
We shall present a
simple numerical example in the final section.  
\biglf
The Direct Meshless Local Petrov Galerkin methods of 
\RScite{mirzaei-schaback:2011-1,mirzaei-et-al:2012-1}
are weak cases of this approach, with other functionals than \eref{eqDuau}
being directly approximated in
terms of values at nodes, and with integrations involving $\bB$ matrices again,
like in all other weak methods. We leave this for a future paper
dealing with all variants of Atluri's Meshless Local Petrov Galerkin technique,
and comparing them to FEM and kernel methods, optimal or not, sparse or not.
\section{Numerical Examples}\RSlabel{SecNumXSCFExa}
To avoid overloading the paper, we present a simple series of examples. 
They all work on the unit circle for simplicity, and in order to include 
finite elements, we have to use triangulations, even if they are not needed for
meshless methods. We start with the standard
discretization of the unit circle into 8 triangles meeting at the origin, 
followed by three standard finite--element refinement
steps halving the edges. The problem \eref{eqLufBug0}
is posed with $L=-\Delta$, and  Dirichlet boundary values are always provided in
the boundary nodes, which are the $y_k$ in \eref{eqdata}.
In all cases, the non--boundary vertices of the triangles  are the nodes where we want to know
the solution, but in the sense of \eref{eqexrs} and for simplicity we only evaluate the recovery
error at the origin $x=0$. 
\biglf
As the FEM variations show, one can work with $f$ values either in barycenters
of triangles or in vertices of the triangulation. We shall evaluate all examples
in both situations, denoting the methods by either *Bary or *Node. Details on
the discretizations are in Table \RSref{TabDisk}:
\begin{itemize}
\item[] $n$: number of Dirichlet boundary data points for $g$ values, 
\item[] $m_{Bary}$: number of triangles and barycentric data points for $f$ values,
\item[] $m_{Node}$: number of vertices and vertex data points for $f$ values,
  including the $n$ Dirichlet boundary vertices,
\item[] DOF: degrees of freedom $=$ number of unknowns $=m_{Node}-n$,
\item[] $h$:  {\em fill
  distance} in the sense of kernel discretizations, describing the maximal
distance of an arbitrary
point of the domain to one of the vertices.
\end{itemize}  
\begin{table}[htp]\centering
\begin{tabular}{||c||r|r|r|r|r||}\hline
Case & $n$ & $m_{Bary}$ &  $m_{Node}$ &DOF  & $h$\\
\hline 
C0 &  8 &   8 &   9 & 1   & 0.2706\\ 
C1 & 16 &  32 &  25 & 9   & 0.1515\\ 
C2 & 32 & 128 &  81 & 49  & 0.0768\\ 
C3 & 64 & 512 & 289 & 225 & 0.0389\\
C4 & 128 & 2048 & 1089 & 961 & 0.0197\\
\hline 
\end{tabular} 
\caption{Discretization data for the examples on the unit disk\RSlabel{TabDisk}}
\end{table}

The plots in Figures \RSref{figOrder3} to \RSref{figOrder7} 
show the errors in Sobolev space of order 3 to 7, with order 3 being
somewhat out of bounds, and the data for the fine discretization C4 being 
polluted by ill--conditioning in various cases. In each figure, the Sobolev space
for error evaluation is fixed, but
the methods might use other kernels. The methods are
\begin{itemize}
\item[] FEMBary: FEM with $f$ data in barycenters
\item[] FEMNode: FEM with $f$ data in nodes
\item[] KansaBary: Unsymmetric collocation  with $f$ data in barycenters, using
  the order 7 Sobolev kernel 
\item[] KansaNode: same with node data
\item[] HOBary: symmetric high-order collocation  with $f$ data in barycenters, using
  the order 7 Sobolev kernel 
\item[] HONode: same with data in nodes. These two coincide with the optimal
  methods, if evaluated on Sobolev space of order 7, see Figure \RSref{figOrder7}. 
\item[] OptBary: optimal method in Sobolev space used for error evaluation, with $f$ data in barycenters
\item[] OptNode: same with data in nodes
\item[] LocNode: a bandwidth 15 method like in section \RSref{SecGDM}, using the
  order 7 Sobolev kernel, data in nodes.
\end{itemize} 
There are serious instabilities that need explanation. They could 
have been avoided in all cases by choosing a smaller kernel scale, but this
would have changed the Sobolev space norm and concealed the 
instabilities. Consequently, the 
kernel scale for Sobolev
space error  evaluation was fixed at 1.0 throughout. This does not seriously
affect the error evaluation, because the latter just consists of a calculation of a
quadratic form. 
\biglf
But it affects the calculation of recovery formulas by kernel
methods, including the optimal ones. In particular, increasing the smoothness of
the construction kernel will increase the instability, if no precautions
like preconditioning 
\RScite{beatson-et-al:1999-1,brown-et-al:2005-1}
are taken. We chose the Kansa* and the HO* kernel
methods to work with the kernel of Sobolev space of order 7 at scale 1,
but smaller scales would have been more stable. Users can easily try different
kernels and scales
for construction and then evaluate in a fixed Sobolev space at a fixed scale
to see which construction scale gives best results. \red{Example is missing}
\biglf
The numerical results support quantitatively what is known from 
experience and partly supported by theory. The piecewise linear FEM technique is adapted
to low--regularity problems, and it shows its second--order convergence 
in all applicable Sobolev spaces from order 4 on. It is clearly inferior
error--wise to
the error--optimal method, which is symmetric collocation in Sobolev space,
and the difference gets larger for increasing Sobolev order, because 
the optimal method behaves like a $p$--FEM and increases its convergence rate
automatically with the Sobolev order. 
Unfortunately, it suffers from
severe ill--conditioning if no precautions are taken, and it 
does not allow sparsity.  However, 
it serves as a standard reference to evaluate error performances of all other methods
using the same data.
\biglf
We now list some observations concerning comparisons.
The optimal method HO* for a fixed high Sobolev order 
performs well also for lower Sobolev order. It adapts automatically to lower
regularity. Other competing methods, like Kansa* for unsymmetric collocation,
behave in the same way. If sparsity is enforced by localizing meshless
collocation \RScite{sarler:2007-1,yao-et-al:2010-1, yao-et-al:2011-1,%
yao-et-al:2012-1,vertnik-sarler:2011-1}, 
one gets competitive methods Loc* that share sparsity with FEM techniques,
but also yield high convergence orders depending on the bandwidth chosen.

\section{Conclusion and Outlook}\RSlabel{SecCC}
The paper provides a tool that allows an explicit and fully computational assessment of 
the error behavior of all linear solvers for all linear PDE problems. 
The exact solution can be unknown, and the error is expressed as a factor
of the (unknown) Sobolev norm of the true solution.
This tool should be applied in many more circumstances, e.g. on special and
awkward domains, for more general differential operators including those
of Computational Mechanics, and for many other
linear solvers, e.g.  MLPG methods and boundary--oriented techniques like the
DRM.
Any application--oriented paper can, in principle, apply this technique and thus
provide a strict worst--case error bound in terms of the Sobolev norm of the
true solution. Examples of single cases with known solutions can never be 
completely satisfactory, but they are the usual practice in application papers.
\biglf
There is a method that always realizes the optimal error, but it needs further
work towards numerical stabilization. All other methods should 
be compared to it, and it is an interesting research challenge to see how close one can
come to the optimal method under sparsity restrictions. 
\begin{figure}[hbt]
\begin{center}
\includegraphics[width=10cm,height=10cm]{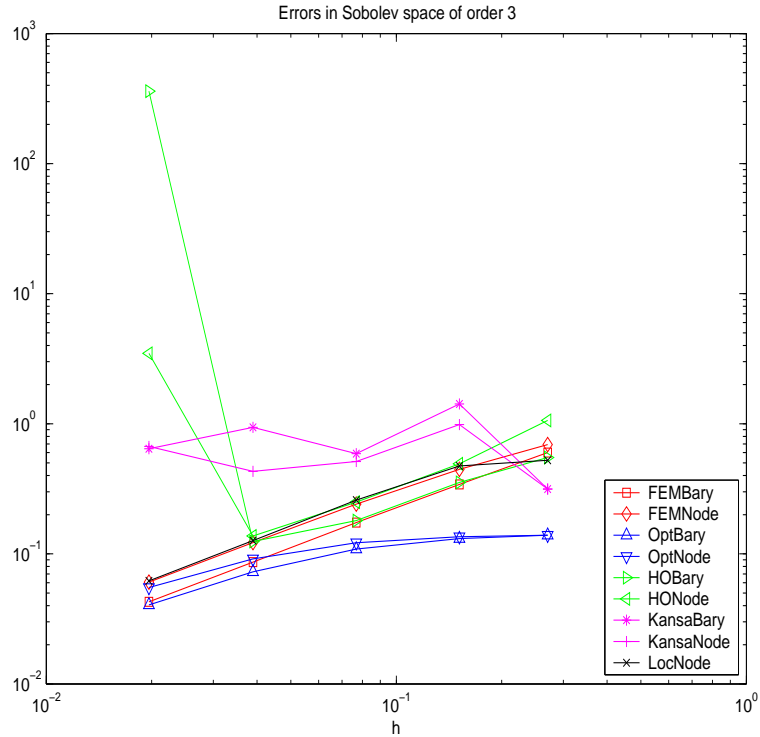}
\caption{Errors for Sobolev order 3}\RSlabel{figOrder3}
\end{center}
\end{figure}

\begin{figure}[hbt]
\begin{center}
\includegraphics[width=10cm,height=10cm]{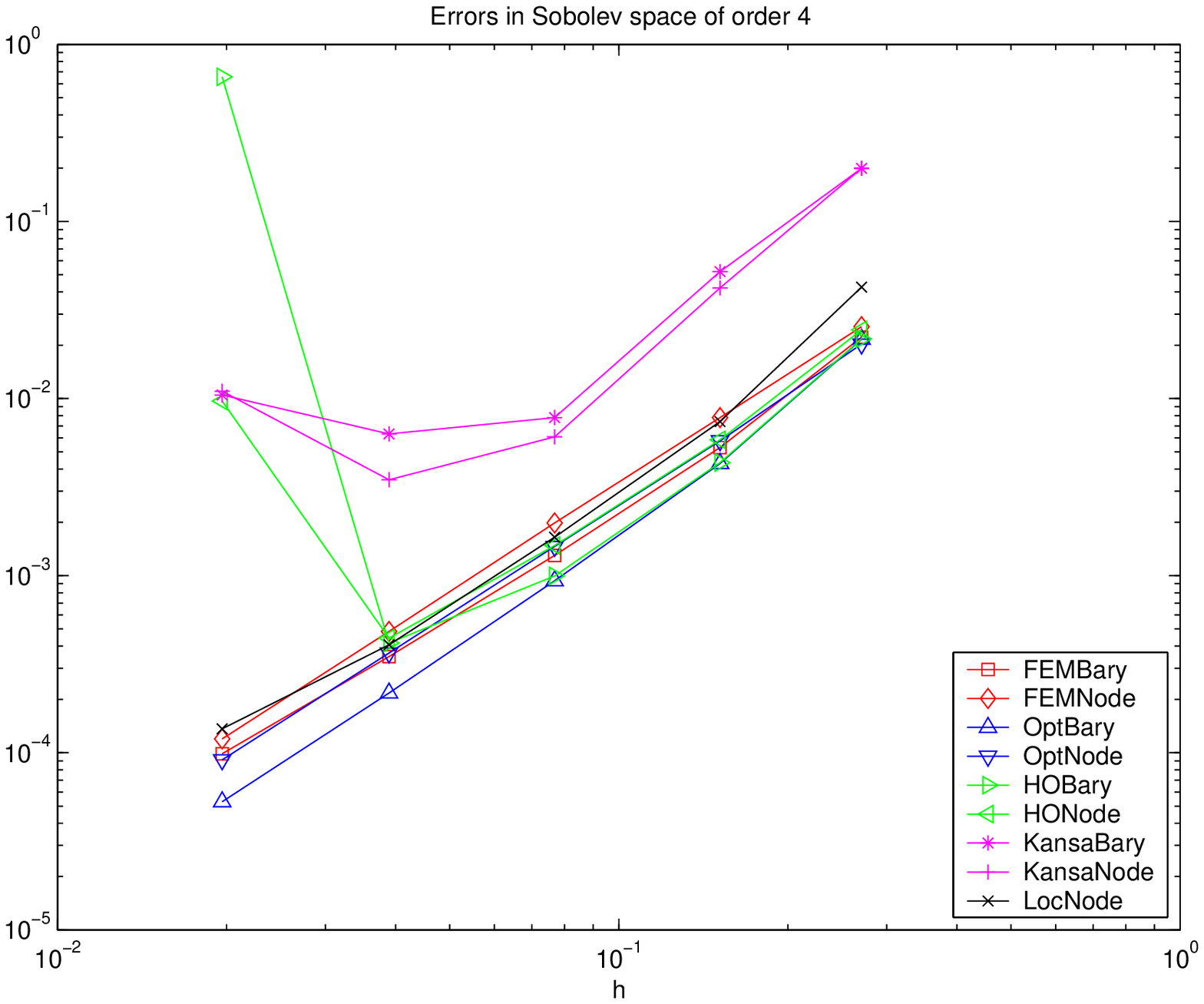}
\caption{Errors for Sobolev order 4}\RSlabel{figOrder4}
\end{center}
\end{figure}

\begin{figure}[hbt]
\begin{center}
\includegraphics[width=10cm,height=10cm]{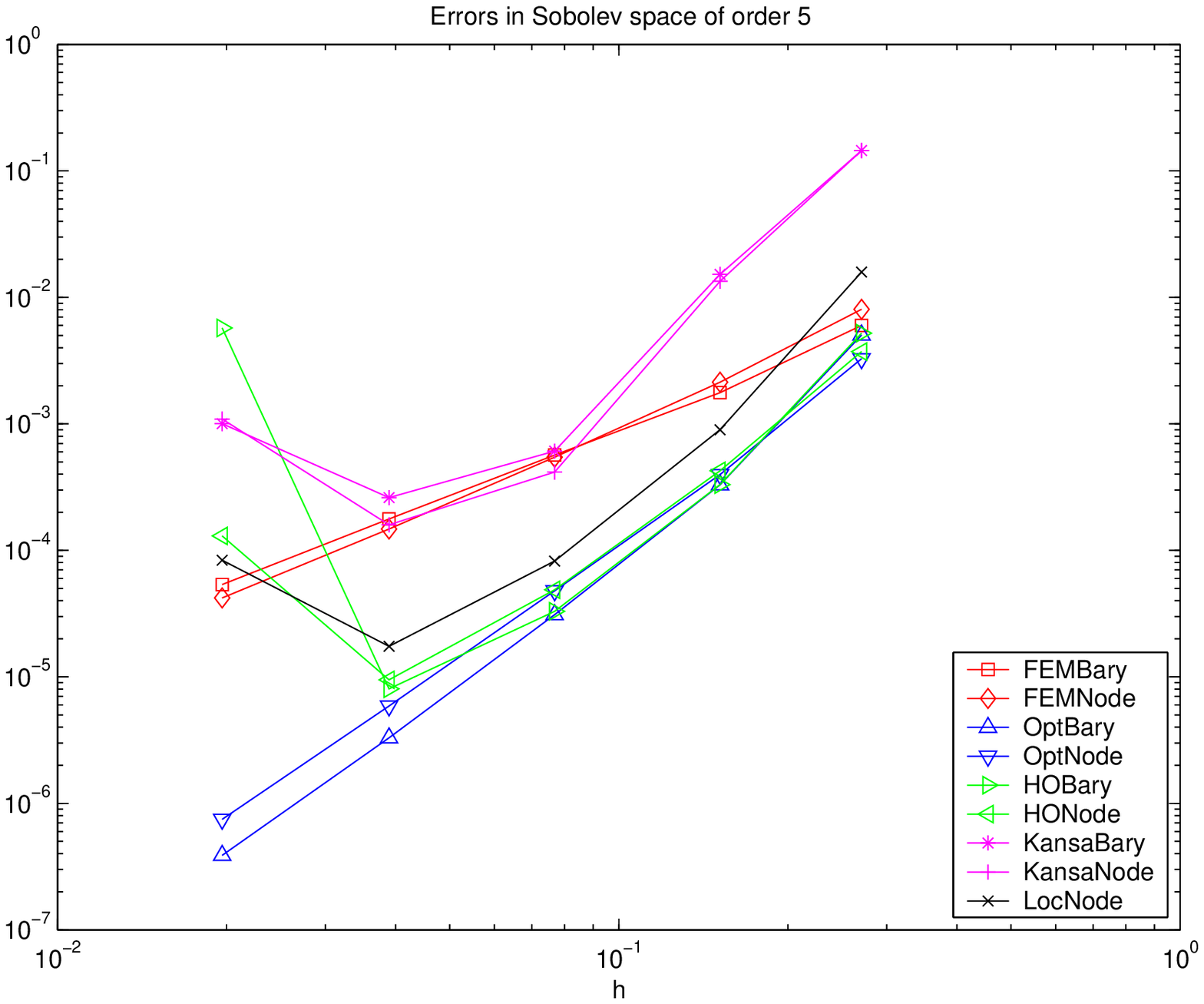}
\caption{Errors for Sobolev order 5}\RSlabel{figOrder5}
\end{center}
\end{figure}

\begin{figure}[hbt]
\begin{center}
\includegraphics[width=10cm,height=10cm]{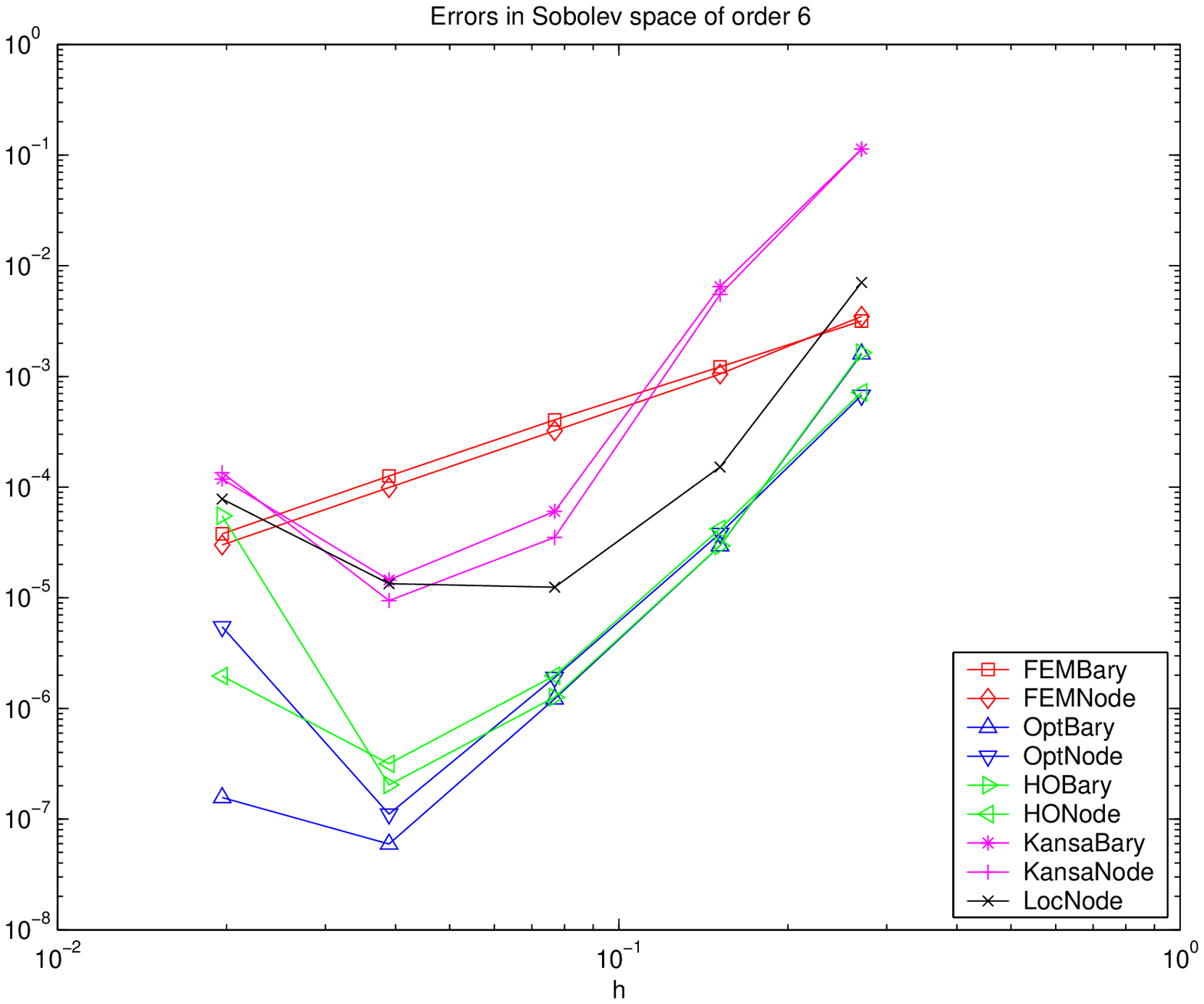}
\caption{Errors for Sobolev order 6}\RSlabel{figOrder6}
\end{center}
\end{figure}

\begin{figure}[hbt]
\begin{center}
\includegraphics[width=10cm,height=10cm]{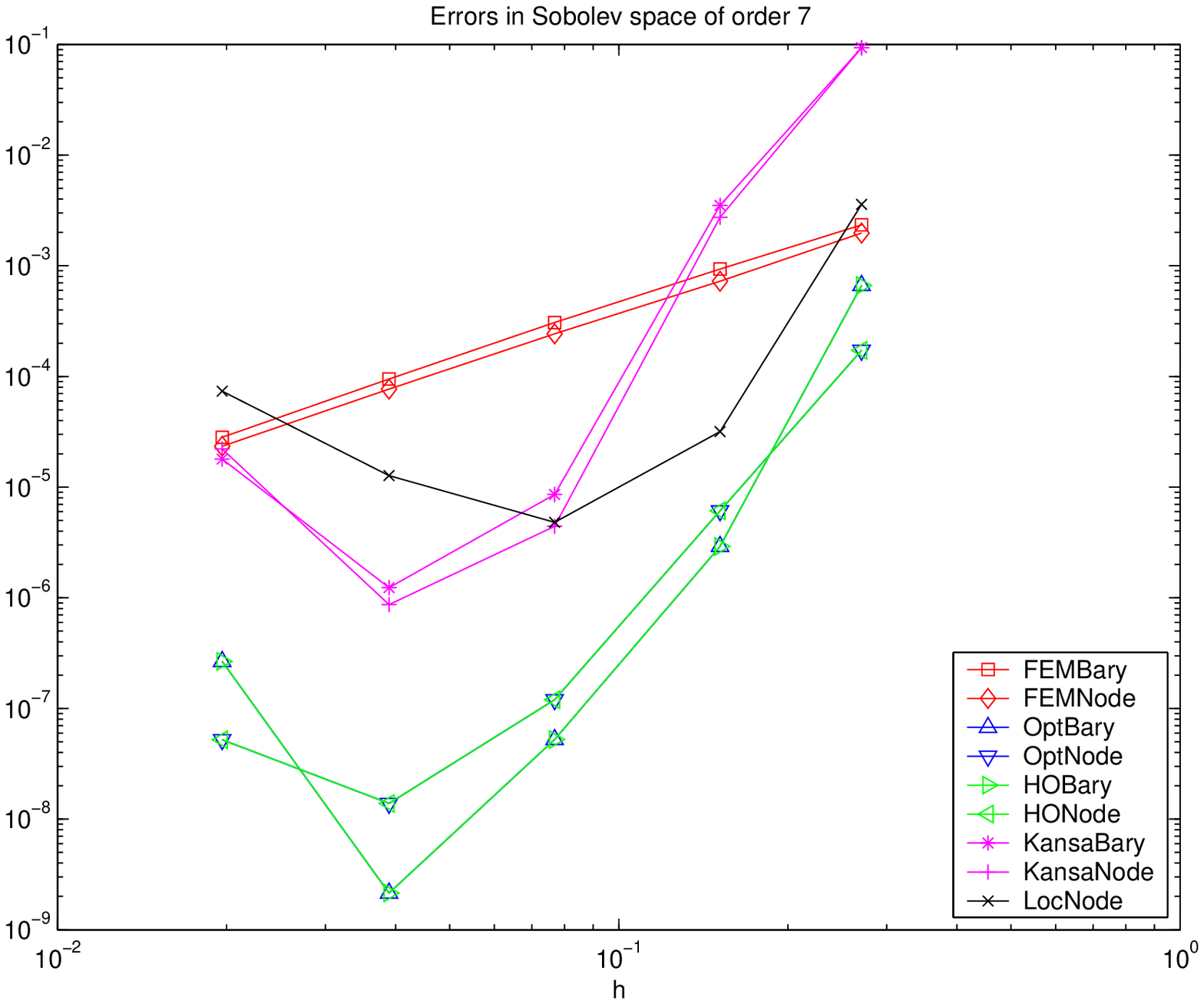}
\caption{Errors for Sobolev order 7}\RSlabel{figOrder7}
\end{center}
\end{figure}

\begin{table}[hbt]
{\small
\begin{tabular}{|r||rrrrr|}
\hline
& C0 & C1 & C2 & C3 & C4 \\
\hline\hline 
h & 2.706e-001 & 1.515e-001 & 7.682e-002 & 3.895e-002 & 1.965e-002 \\ 
FEMBary & 6.011e-001 & 3.415e-001 & 1.737e-001 & 8.635e-002 & 4.272e-002 \\ 
FEMNode & 6.912e-001 & 4.457e-001 & 2.411e-001 & 1.220e-001 & 6.063e-002 \\ 
KansaBary & 3.150e-001 & 1.422e+000 & 5.888e-001 & 9.385e-001 & 6.421e-001 \\ 
KansaNode & 3.150e-001 & 9.847e-001 & 5.131e-001 & 4.303e-001 & 6.696e-001 \\  
HOBary & 5.510e-001 & 3.538e-001 & 1.802e-001 & 1.245e-001 & 3.605e+002 \\  
HONode & 1.060e+000 & 4.918e-001 & 2.498e-001 & 1.371e-001 & 3.475e+000 \\  
OptBary & 1.390e-001 & 1.309e-001 & 1.090e-001 & 7.280e-002 & 4.051e-002  \\ 
OptNode & 1.391e-001 & 1.353e-001 & 1.218e-001 & 9.145e-002 & 5.508e-002 \\ 
LocNode & 5.236e-001 & 4.734e-001 & 2.591e-001 & 1.264e-001 & 6.196e-002  \\ 
\hline
\end{tabular}
}
\caption{Raw data for Sobolev order 3\RSlabel{TabSob3}}
\end{table}
\begin{table}[hbt]
{\small
\begin{tabular}{|r||rrrrr|}
\hline
& C0 & C1 & C2 & C3 & C4 \\
\hline\hline 
h & 2.706e-001 & 1.515e-001 & 7.682e-002 & 3.895e-002 & 1.965e-002 \\ 
FEMBary & 2.223e-002 & 5.296e-003 & 1.301e-003 & 3.494e-004 & 9.856e-005 \\ 
FEMNode & 2.547e-002 & 7.813e-003 & 1.984e-003 & 4.844e-004 & 1.199e-004 \\ 
KansaBary & 1.995e-001 & 5.219e-002 & 7.806e-003 & 6.309e-003 & 1.049e-002 \\ 
KansaNode & 1.995e-001 & 4.215e-002 & 6.083e-003 & 3.478e-003 & 1.101e-002 \\  
HOBary & 2.173e-002 & 4.353e-003 & 9.958e-004 & 4.171e-004 & 6.565e-001 \\  
HONode & 2.446e-002 & 5.855e-003 & 1.479e-003 & 4.436e-004 & 9.688e-003 \\  
OptBary & 2.163e-002 & 4.308e-003 & 9.355e-004 & 2.174e-004 & 5.287e-005  \\ 
OptNode & 2.029e-002 & 5.770e-003 & 1.463e-003 & 3.653e-004 & 9.117e-005 \\ 
LocNode & 4.265e-002 & 7.420e-003 & 1.651e-003 & 4.069e-004 & 1.364e-004  \\ 
\hline
\end{tabular}
}
\caption{Raw data for Sobolev order 4\RSlabel{TabSob4}}
\end{table}
\begin{table}[hbt]
{\small
\begin{tabular}{|r||rrrrr|}
\hline
& C0 & C1 & C2 & C3 & C4 \\ 
\hline\hline 
h & 2.706e-001 & 1.515e-001 & 7.682e-002 & 3.895e-002 & 1.965e-002 \\ 
FEMBary & 6.010e-003 & 1.762e-003 & 5.658e-004 & 1.775e-004 & 5.376e-005 \\ 
FEMNode & 8.042e-003 & 2.134e-003 & 5.451e-004 & 1.477e-004 & 4.202e-005 \\ 
KansaBary & 1.447e-001 & 1.522e-002 & 6.089e-004 & 2.612e-004 & 1.005e-003 \\ 
KansaNode & 1.447e-001 & 1.345e-002 & 4.163e-004 & 1.590e-004 & 1.092e-003 \\  
HOBary & 5.223e-003 & 3.307e-004 & 3.299e-005 & 8.055e-006 & 5.733e-003 \\  
HONode & 3.723e-003 & 4.258e-004 & 4.881e-005 & 9.465e-006 & 1.303e-004 \\  
OptBary & 5.031e-003 & 3.304e-004 & 3.106e-005 & 3.298e-006 & 3.873e-007  \\ 
OptNode & 3.278e-003 & 3.977e-004 & 4.798e-005 & 5.890e-006 & 7.463e-007 \\ 
LocNode & 1.591e-002 & 8.987e-004 & 8.243e-005 & 1.742e-005 & 8.366e-005  \\ 
\hline
\end{tabular}
}
\caption{Raw data for Sobolev order 5\RSlabel{TabSob5}}
\end{table}
\begin{table}[hbt]
{\small
\begin{tabular}{|r||rrrrr|}
\hline
& C0 & C1 & C2 & C3 & C4 \\ 
\hline\hline
h & 2.706e-001 & 1.515e-001 & 7.682e-002 & 3.895e-002 & 1.965e-002 \\ 
FEMBary & 3.167e-003 & 1.219e-003 & 4.051e-004 & 1.262e-004 & 3.780e-005 \\ 
FEMNode & 3.485e-003 & 1.054e-003 & 3.226e-004 & 9.917e-005 & 3.003e-005 \\ 
KansaBary & 1.135e-001 & 6.470e-003 & 6.036e-005 & 1.451e-005 & 1.186e-004 \\ 
KansaNode & 1.135e-001 & 5.501e-003 & 3.513e-005 & 9.405e-006 & 1.349e-004 \\  
HOBary & 1.641e-003 & 2.955e-005 & 1.262e-006 & 2.034e-007 & 5.496e-005 \\  
HONode & 7.159e-004 & 4.206e-005 & 1.983e-006 & 3.158e-007 & 1.965e-006 \\  
OptBary & 1.600e-003 & 2.949e-005 & 1.222e-006 & 5.964e-008 & 1.566e-007  \\ 
OptNode & 6.730e-004 & 3.811e-005 & 1.901e-006 & 1.110e-007 & 5.458e-006 \\ 
LocNode & 7.056e-003 & 1.521e-004 & 1.247e-005 & 1.342e-005 & 7.800e-005  \\ 
\hline
\end{tabular}
}
\caption{Raw data for Sobolev order 6\RSlabel{TabSob6}}
\end{table}
\begin{table}[hbt]
{\small
\begin{tabular}{|r||rrrrr|}
\hline
& C0 & C1 & C2 & C3 & C4 \\ 
\hline\hline
h & 2.706e-001 & 1.515e-001 & 7.682e-002 & 3.895e-002 & 1.965e-002 \\ 
FEMBary & 2.337e-003 & 9.327e-004 & 3.075e-004 & 9.490e-005 & 2.820e-005 \\ 
FEMNode & 1.975e-003 & 7.244e-004 & 2.425e-004 & 7.684e-005 & 2.339e-005 \\ 
KansaBary & 9.332e-002 & 3.506e-003 & 8.576e-006 & 1.232e-006 & 1.792e-005 \\ 
KansaNode & 9.332e-002 & 2.743e-003 & 4.428e-006 & 8.672e-007 & 2.207e-005 \\  
HOBary & 6.633e-004 & 2.916e-006 & 5.232e-008 & 2.147e-009 & 2.658e-007 \\  
HONode & 1.722e-004 & 6.116e-006 & 1.202e-007 & 1.390e-008 & 5.223e-008 \\  
OptBary & 6.633e-004 & 2.916e-006 & 5.232e-008 & 2.147e-009 & 2.658e-007  \\ 
OptNode & 1.722e-004 & 6.116e-006 & 1.202e-007 & 1.390e-008 & 5.224e-008 \\ 
LocNode & 3.590e-003 & 3.179e-005 & 4.802e-006 & 1.273e-005 & 7.360e-005  \\ 
\hline
\end{tabular}
}
\caption{Raw data for Sobolev order 7\RSlabel{TabSob7}}
\end{table}

\bibliographystyle{plain}

\end{document}